\documentclass[a4paper,11pt]{amsart}
\usepackage{graphicx}
\usepackage{mathptmx}
\usepackage{amsmath}
\usepackage{amssymb}
\usepackage{enumitem}
\usepackage{xcolor}

\newmuskip\pFqmuskip

\newcommand*\pFq[6][8]{%
  \begingroup 
  \pFqmuskip=#1mu\relax
  \mathcode`=\string"8000
  \begingroup\lccode`\~=`\,
  \lowercase{\endgroup\let~}\pFqcomma
  F^{#2}_{#3}{\left(\genfrac..{0pt}{}{#4}{#5}\bigg|#6\right)}%
  \endgroup
}
\newcommand{\pFqcomma}{\mskip\pFqmuskip}

\newtheorem{theorem}{Theorem}

\newtheorem{definition}[theorem]{Definition}
\newtheorem{remark}[theorem]{Remark}

\begin{document}

\title[]{Study on discrete degenerate Bell distributions with two parameters}

\author{Taekyun Kim $^{1}$}
\address{Department of Mathematics, Kwangwoon University, Seoul 139-701, Republic of Korea}
\email{tkkim@kw.ac.kr}

\author{Dae San Kim $^{2}$}
\address{Department of Mathematics, Sogang University, Seoul 121-742, Republic of Korea}
\email{dskim@sogang.ac.kr}

\author{Hye Kyung Kim $^{3,*}$}
\address{Department of Mathematics Education, Daegu Catholic University, Gyeongsan 38430, Republic of Korea}
\email{hkkim@cu.ac.kr}

\thanks{* is corresponding author}

\begin{abstract}
Recently, Freud-Rodriguez proposed a new counting process which is called the Bell-Touchard process and based on the Bell-Touchard probability distribution. This process was developed to solve the problem of rare events hypothesis which is one of the limitations of the Poisson process. In this paper, we consider the discrete degenerate Bell distributions and the degenerate Bell process which are 'degenerate versions' of the Bell-Touchard probability distributions and the Bell-Touchard process, respectively. We investigate several properties of the degenerate Bell distribution. We introduce the degenerate Bell process by giving two equivalent definitions and show one method of constructing a new infinite family of degenerate Bell process out of a given infinite family of degenerate Bell process.
\end{abstract}

 \maketitle

\markboth{\centerline{\scriptsize Study on discrete degenerate Bell distributions with two parameters}}
{\centerline{\scriptsize  T. Kim, D. S. Kim and H. K. Kim}}


\section{Introduction}
We have witnessed in recent years that studying various degenerate versions of some special polynomials and numbers yields many interesting and fruitful results. This explorations for degenerate versions started from the pioneering work of Carlitz on the degenerate Bernoulli and degenerate Euler polynomials. \par
It turns out that the Bell-Touchard distribution with parameters $(\alpha, \theta)$ can be used quite effectively for modeling count data (see \cite{2}).
In \cite{5}, Freud and Rodriguez defined the Bell-Touchard process which is based the Bell-Touchard probability distribution. It is introduced to contribute the formulation of mathematical models where the rare events hypothesis is not suitable. The aim of this paper is to study the discrete degenerate Bell distributions and the degenerate Bell process which are degenerate versions of the Bell-Touchard probability distributions and the Bell-Touchard process, respectively. We investigate several properties of the degenerate Bell-Touchard probability distributions. Then we give two equivalent definitions of the degenerate Bell-Touchard process. We show that the infinite family of degenerate Bell process formed by the partial sums of an infinite family of degenerate Bell process with the same $\theta$ parameter is of the same nature.\par
In more detail, the outline of this paper is as follows. In Section 1, we recall the degenerate exponential functions, the degenerate Stirling numbers of the second kind and the degenerate Bell polynomials. Then we remind the reader of the counting process, independent increments and stationary increments. We also recall the Poisson process. Section 2 is the main result of this paper. We derive the Dobinski-like formula for the degenerate Bell polynomials. We introduce the degenerate Bell distributions with two parameters $(\alpha,\theta)$. Then we find the probability generating function of such distributions in Theorem 1. We deduce the parameters for the sum of a finite number of independent degenerate Bell distributions with the same $\theta$ parameter in Theorem 2. We determine the moment generating function of the degenerate Bell distributions in Theorem 3. In Theorem 4, we compute the expectation and the variance of the degenerate Bell distributions. We define the degenerate Bell process with parameters $(\alpha,\theta)$. We introduce the degenerate Bell process $\{N_{\lambda}(t)|\,t \ge 0 \}$ with parameters $(\alpha,\theta)$ as s counting process satisfying conditions (i), (ii) and (iii). Here the condition (iii) specifies only the linear term of $P\{N_{\lambda}(t)=k\}$, the probability of $k$ events occurring in the interval $[0,t)$. However, we show that the conditions (i), (ii) and (iii) together determine $P\{N_{\lambda}(t)=k\}$ completely in Theorem 5, so that we come up with another (equivalent) definition for the degenerate Bell process in Definition 6. In Theorem 7, we show that, for a given infinite family of degenerate Bell process with the same $\theta$ parameter, the infinite family of degenerate Bell process formed by their partial sums are also of the same nature. Finally, we conclude our paper in Section 3. In the rest of this section, we recall the necessary facts that are needed throughout this paper. \par

\vspace{0.1in}
For any nonzero $\lambda \in \mathbb{R}$, the degenerate exponentials are defined by
\begin{equation}\label{eq01}
\begin{split}
&e_\lambda^x(t)=(1+\lambda t)^{\frac{x}{\lambda}}=\sum_{n=0}^\infty \frac{(x)_{n,\lambda}}{n!}t^n, \ e_\lambda(t)=e_\lambda^1(t), \quad(\text{see [7-11]}),
\end{split}
\end{equation}
where
$(x)_{0,\lambda}=1, \  (x)_{n,\lambda}=x(x-\lambda) \cdots (x-(n-1)\lambda), \quad (n\geq 1).$ \\
Note that $\lim_{\lambda\rightarrow 0}(x)_{n,\lambda}=x^n$, and $\lim_{\lambda\rightarrow 1}(x)_{n,\lambda}=(x)_n$,
\noindent where $(x)_0=1, (x)_n=x(x-1)\cdots(x-n+1), \quad (n\geq 1)$.


The Stirling numbers of the second kind are defined by 
\begin{equation}\label{eq03}
\begin{split}
\frac{1}{k!}(e^t-1)^k=\sum_{n=k}^\infty S_2(n,k)\frac{t^n}{n!}, \quad {(\text {see \cite{3,11,8}})}.
\end{split}
\end{equation}
It is well known that the Bell polynomials are defined by 
\begin{equation}\label{eq02}
\begin{split}
e^{x(e^t-1)}=\sum_{n=0}^\infty \phi_n(x)\frac{t^n}{n!}, \quad {(\text {see \ [2-5]})}.
\end{split}
\end{equation}
From \eqref{eq03} and \eqref{eq02}, we note that
\begin{equation}\label{eq04}
\begin{split}
\phi_{n}(x)=\sum_{k=0}^n S_2(n,k)x^k, \quad (n\geq 0), \quad {(\text {see \cite{9,10}})}.
\end{split}
\end{equation}\par
The degenerate Stirling numbers of the second kind are defined by
\begin{equation}\label{eq05}
\begin{split}
(x)_{n,\lambda}=\sum_{k=0}^n S_{2,\lambda}(n,k)(x)_k, \quad (n\geq 0), \quad {(\text {see \ [7-11]})}.
\end{split}
\end{equation}
From \eqref{eq05}, we note that
\begin{equation}\label{eq06}
\begin{split}
\frac{1}{k!}(e_\lambda(t)-1)^k=\sum_{n=k}^\infty S_{2,\lambda}(n,k)\frac{t^n}{n!}, \quad {(\text {see \cite{10}})}.
\end{split}
\end{equation}
Note that $\lim_{\lambda\rightarrow0}S_{2,\lambda}(n,k)=S_2(n,k), \quad (n,k\geq0)$. \par
Recently, the degenerate Bell polynomials are defined by Kim-Kim as
\begin{equation}\label{eq07}
\begin{split}
e^{x(e_\lambda(t)-1)}=\sum_{n=0}^\infty \phi_{n,\lambda}(x)\frac{t^n}{n!}, \quad {(\text {see \cite{10}})}.
\end{split}
\end{equation}
Note that the degenerate Bell polynomials satisfy the binomial identity:
\begin{equation}\label{eq07-1}
\begin{split}
\phi_{n,\lambda}(x+y)=\sum_{k=0}^{n}\binom{n}{k}\phi_{k,\lambda}(x)\phi_{n-k,\lambda}(y).
\end{split}
\end{equation}
By \eqref{eq06} and \eqref{eq07}, we get
\begin{equation}\label{eq08}
\begin{split}
\phi_{n,\lambda}(x)=\sum_{k=0}^n S_{2,\lambda}(n,k)x^k, \quad {(\text {see \,[7-11]})}.
\end{split}
\end{equation} \par

A stochastic process $\{N(t)|\,t\geq0\}$ is said to be a {\it{counting process}} if $N(t)$ represents the total number of `events' that occur by time $t$. \\
From its definition, we note that for a counting process $N(t)$ must satisfy {(\text {see \cite{1,2,5,6,12}})}:

\begin{itemize}
\item[(i)] $N(t)\geq0$, \\
\item[(ii)] $N(t)$ is integer valued, \\
\item[(iii)] if $s<t$, then $N(s)\leq N(t)$, \\
\item[(iv)] for $s<t$, $N(t)- N(s)$ equals the number of events that occur in the interval $(s, \ t].$
\end{itemize}

A counting process is said to {\it{ possess independent increments}} if the numbers of events that occur in disjoint time intervals are independent.\\
 A counting process is said to {\it{possess stationary increments }}if the distribution of the number of events that occur in any interval of time depends only on the length of the time interval.

The counting process $\{N(t)|\,t\geq0\}$ is said to be a {\it{Poisson process having rate $\alpha \ (\alpha >0)$}} if
\begin{itemize}
\item[(i)] $N(0)=0$, \\
\item[(ii)] the process has independent increments, \\
\item[(iii)] the number of events in any interval of length $t$ is Poisson distributed with mean $\alpha t$, \quad {(\text {see \cite{1}})}.
\end{itemize}

That is, for all $s,\ t\geq 0$,
\begin{equation*}
\begin{split}
P\{N(t+s)-N(s)=n\}=e^{-\alpha t}\frac{(\alpha t)^n}{n!}.
\end{split}
\end{equation*}

In this paper, we introduce a new degenerate Bell discrete distribution with two parameters and propose a new counting process based on the degenerate Bell probability distribution, naming it the degenerate Bell process.

\medskip

\section{Discrete degenerate Bell distributions with two parameters}
In this section, we assume that $\lambda\in (0,1]$.
From \eqref{eq07}, we note that
\begin{equation}\label{eq09}
\begin{split}
\sum_{n=0}^\infty \phi_{n,\lambda}(x)\frac{t^n}{n!}&=e^{x(e_\lambda(t)-1)}=e^{-x}e^{xe_\lambda(t)} \\
&=e^{-x}\sum_{k=0}^\infty \frac{x^k}{k!}e_\lambda^k(t)=\sum_{n=0}^\infty \bigg(e^{-x}\sum_{k=0}^\infty\frac{(k)_{n,\lambda}}{k!}x^{k}\bigg)\frac{t^n}{n!}.
\end{split}
\end{equation}
Comparing the coefficients on both sides of \eqref{eq09}, we get the Dobinski-like formula:
\begin{equation}\label{eq10}
\begin{split}
\phi_{n,\lambda}(x)=e^{-x}\sum_{k=0}^\infty \frac{(k)_{n,\lambda}}{k!}x^{k}, \quad (n\geq0), \quad {(\text {see \cite{10}})}.
\end{split}
\end{equation}
In particular, for $x=1$, we have
\begin{equation}\label{eq11}
\begin{split}
\phi_{n,\lambda}=\phi_{n,\lambda}(1)=\frac{1}{e}\sum_{k=0}^\infty \frac{(k)_{n,\lambda}}{k!}, \quad (n\geq0),
\end{split}
\end{equation}
which are called the degenerate Bell numbers.

Now, we consider the degenerate Bell random variable with parameters $\alpha$ and $\theta$.\\
A discrete random variable $X$ has a {\it{degenerate Bell distribution}} with parameters $(\alpha, \theta)\in \mathbb{R}^2$ if its probability mass function is given by
\begin{equation}\label{eq12}
\begin{split}
p(k)=P\{X=k\}=e^{-\alpha(e_\lambda(\theta)-1)}\frac{\theta^k}{k!}\phi_{k,\lambda}(\alpha), \quad (k\geq0),
\end{split}
\end{equation}
which is denoted by $X\sim DB_\lambda(\alpha,\theta)$.
Note that
\begin{equation*}
\begin{split}
\sum_{k=0}^\infty p(k)=\sum_{k=0}^\infty P\{X=k\}&=e^{-\alpha(e_\lambda(\theta)-1)} \sum_{k=0}^\infty \frac{\phi_{k,\lambda(\alpha)}}{k!}\theta^k \\
&=e^{-\alpha(e_\lambda(\theta)-1)}e^{\alpha(e_\lambda(\theta)-1)}=1.
\end{split}
\end{equation*}

Let $X\sim DB_\lambda(\alpha,\theta)$. Then the probability generating function $G_{X}(t)$ of $X$ is given by
\begin{equation}\label{eq13}
\begin{split}
G_{X}(t)&=E[t^X]=\sum_{n=0}^\infty P\{X=n\}t^n \\
&=e^{-\alpha(e_\lambda(\theta)-1)}\sum_{n=0}^\infty\frac{(\theta t)^n}{n!}\phi_{n,\lambda}(\alpha) \\
&=e^{-\alpha(e_\lambda(\theta)-1)}e^{\alpha(e_\lambda(\theta t)-1)} =e^{\alpha(e_\lambda(\theta t)-e_\lambda(\theta))}.
\end{split}
\end{equation}

Therefore, by \eqref{eq13}, we obtain the following theorem.

\begin{theorem}
Let $X\sim DB_\lambda(\alpha,\theta)$. Then the probability generating function $G_{X}(t)$ of $X$ is given by
\begin{equation*}
\begin{split}
G_X(t)=e^{\alpha(e_\lambda(\theta t)-e_\lambda(\theta))}.
\end{split}
\end{equation*}
\end{theorem}

\medskip

Let $\{X_i\}_{i=1}^n$ be the sequence of independent random variables with $X_i\sim DB_\lambda(\alpha_i,\theta)$,\\
 and let $Y=\sum_{i=1}^nX_i$. Then we have
\begin{equation}\label{eq14}
\begin{split}
G_Y(t)&=E[t^Y] =E[t^{\sum_{i=1}^nX_i}]=\Pi_{i=1}^nE[t^{X_i}] \\
&=\Pi_{i=1}^n e^{\alpha_i(e_\lambda(\theta t)-e_\lambda(\theta))} 
=e^{\sum_{i=1}^n\alpha_i(e_\lambda(\theta t)-e_\lambda(\theta))}.
\end{split}
\end{equation}

Therefore, by \eqref{eq13} and  \eqref{eq14}, we obtain the following theorem.

\begin{theorem}
Let $\{X_i\}_{i=1}^n$ be the sequence of independent random variables with $X_i\sim DB_\lambda(\alpha_i,\theta)$.
Then we have
\begin{equation*}
\begin{split}
\sum_{i=1}^nX_i\sim DB_\lambda\Big(\sum_{i=1}^n \alpha_i,\theta\Big).
\end{split}
\end{equation*}
\end{theorem}

\medskip

Let $X$ be the discrete random variable with probability mass function given by $p(k)=P\{X=k\}$. \\
Then the moments of $X$ are defined by
\begin{equation*}
\begin{split}
E[X^n]=\sum_{k=0}^\infty k^np(k)=\sum_{k=0}^\infty k^nP\{X=k\}, \quad (n\geq0).
\end{split}
\end{equation*}

The moment generating function of $X$ is given by
\begin{equation*}
\begin{split}
F_X(t)=E[e^{Xt}]=\sum_{n=0}^\infty E[X^n]\frac{t^n}{n!}=\sum_{n=0}^\infty e^{nt}P_r\{X=n\} , \quad {(\text {see \cite{12,13}})}.
\end{split}
\end{equation*}

Let $X \sim DB_\lambda(\alpha, \theta)$. Then, by \eqref{eq07} and \eqref{eq12}, the moment generating function of $X$ is given by
\begin{equation}\label{eq15}
\begin{split}
F_X(t)=E[e^{Xt}]&=\sum_{n=0}^\infty e^{nt}P\{X=n\} \\
&=e^{-\alpha(e_\lambda(\theta)-1)}\sum_{n=0}^\infty e^{nt}\frac{\theta^n \phi_{n,\lambda}(\alpha)}{n!} \\
&=e^{-\alpha(e_\lambda(\theta)-1)}e^{\alpha(e_\lambda(e^t\theta)-1)} =e^{\alpha(e_\lambda(e^t\theta)-e_\lambda(\theta))}.
\end{split}
\end{equation}

Therefore, by \eqref{eq15}, we obtain the following theorem.
\begin{theorem}
For $X \sim DB_\lambda(\alpha,\theta)$, let $F_X(t)=\sum_{n=0}^\infty E(X^n)\frac{t^n}{n!}$ be the moment generating function of $X$. Then we have
\begin{equation*}
\begin{split}
F_X(t)=e^{\alpha(e_\lambda(e^t\theta)-e_\lambda(\theta))}.
\end{split}
\end{equation*}
\end{theorem}

\medskip

Let $X\sim DB_{\lambda}(\alpha, \theta)$. Then we have
\begin{equation}\label{eq26}
\begin{split}
E[X]=\sum_{k=0}^{\infty}k P\{X=k\}&=\sum_{k=0}^\infty k\cdot e^{-\alpha(e_{\lambda}(\theta)-1)}\frac{\theta^k}{k!}\phi_{k,\lambda}(\alpha)\\
&=e^{-\alpha(e_{\lambda}(\theta)-1)}\theta \frac{\partial}{\partial\theta}\sum_{k=0}^\infty\frac{\theta^k}{k!}\phi_{k,\lambda}(\alpha)\\
&=e^{-\alpha(e_{\lambda}(\theta)-1)}\theta \frac{\partial}{\partial\theta}\Big(e^{\alpha(e_{\lambda}(\theta)-1)}\Big)\\
&=e^{-\alpha(e_{\lambda}(\theta)-1)}\theta\alpha e_{\lambda}^{1-\lambda}(\theta)(e^{\alpha(e_{\lambda}(\theta)-1)})=\theta\alpha e_{\lambda}^{1-\lambda}(\theta),
\end{split}
\end{equation}

and
\begin{equation}\label{eq27}
\begin{split}
E[X^2]&=\sum_{k=0}^\infty k^2e^{-\alpha(e_{\lambda}(\theta)-1)}\frac{\theta^k}{k!}\phi_{k,\lambda}(\alpha)\\
&=e^{-\alpha(e_{\lambda}(\theta)-1)}\bigg(\theta\frac{\partial}{\partial\theta}\bigg)^2\sum_{k=0}^\infty\frac{\theta^k}{k!}\phi_{k,\lambda}(\alpha)\\
&=e^{-\alpha(e_{\lambda}(\theta)-1)}\theta\frac{\partial}{\partial\theta}\bigg(\theta\frac{\partial}{\partial\theta}e^{\alpha(e_{\lambda}(\theta)-1)}\bigg)\\
&=e^{-\alpha(e_{\lambda}(\theta)-1)}\theta\frac{\partial}{\partial\theta}\Big(\theta\alpha e_{\lambda}^{1-\lambda}(\theta))e^{\alpha(e_{\lambda}(\theta)-1)}\Big)\\
&=e^{-\alpha(e_{\lambda}(\theta)-1)}\theta\alpha \Big\{\Big(e_{\lambda}^{1-\lambda}(\theta)+\theta(1-\lambda)e_{\lambda}^{1-2\lambda}(\theta)\Big)e^{\alpha(e_{\lambda}(\theta)-1)}\\
&\hspace{5cm}+\theta e_{\lambda}^{1-\lambda}(\theta)\alpha e_{\lambda}^{1-\lambda}(\theta)e^{\alpha(e_{\lambda}(\theta)-1)}\Big\}\\
&=\theta\alpha \big(1+\theta(1-\lambda)e_{\lambda}^{-\lambda}(\theta)\big)e_{\lambda}^{1-\lambda}(\theta)+\theta^2\alpha^2e_{\lambda}^{2(1-\lambda)}(\theta).
\end{split}
\end{equation}

From \eqref{eq26} and \eqref{eq27}, we have
\begin{equation}\label{eq28}
\begin{split}
\rm{Var}(X)&=E[x^2]-(E[x])^2\\
&=\theta\alpha(1+\theta(1-\lambda)e_{\lambda}^{-\lambda}(\theta))e_{\lambda}^{1-\lambda}(\theta)+\theta^2\alpha^2e_{\lambda}^{2(1-\lambda)}(\theta)-\theta^2\alpha^2e_{\lambda}^{2(1-\lambda)}(\theta)\\
&=\theta\alpha(1+\theta(1-\lambda)e_{\lambda}^{-\lambda}(\theta))e_{\lambda}^{1-\lambda}(\theta).
\end{split}
\end{equation}

Therefore, by \eqref{eq28}, we obtain the following theorem.
\begin{theorem}
Let $X\sim DB_{\lambda}(\alpha, \theta)$. Then we have
\begin{equation*}
\begin{split}
E[X]=\theta\alpha e_{\lambda}^{1-\lambda}(\theta),\ \ {\rm{and}} \ \ \
{\rm{Var}}(X)=\theta\alpha \big(1+\theta(1-\lambda)e_{\lambda}^{-\lambda}(\theta)\big)e_{\lambda}^{1-\lambda}(\theta).
\end{split}
\end{equation*}
\end{theorem}

\medskip

A counting process $\{N_\lambda(t)|t\geq0\}$ is said to be  a {\it{degenerate Bell process}} with parameters $(\alpha, \theta) \in \mathbb{R}^2_+$ if the following assumption hold:
\begin{itemize}
\item[(i)] $N_\lambda(0)=0$, \\
\item[(ii)] $\{N_\lambda(t)|t\geq0\}$ has stationary and independent increments, \\
\item[(iii)] $P\{N_\lambda(t+s)-N_\lambda(t)=k\}=\alpha s\frac{(1)_{k,\lambda}}{k!}\theta^k+o(s)$, where $k\in\mathbb{N}$ and $s, \ t\geq0$.
\end{itemize}

\medskip

Let $g(t)=E[\exp(-xN_\lambda(t))]$. Then
\begin{equation}\label{eq16}
\begin{split}
g(t+s)&=E[\exp(-xN_\lambda(t+s))] \\
&=E[\exp(-x(N_\lambda(t+s)-N_\lambda(t)+N_\lambda(t))] \\
&=E[\exp(-x(N_\lambda(t))\exp(N_\lambda(t+s)-N_\lambda(t)))] \\
&=E[\exp(-xN_\lambda(t))]E[\exp(-x(N_\lambda(t+s)-N_\lambda(t)))] \\
&=g(t)E[\exp(-xN_\lambda(s))].
\end{split}
\end{equation}

Now, we note that
\begin{equation}\label{eq17}
\begin{split}
\sum_{k=0}^\infty P\{N_\lambda(s)=k\}=1=P\{N_\lambda(s)=0\}+\sum_{k=1}^\infty P\{N_\lambda(s)=k\}.
\end{split}
\end{equation}

Thus, by \eqref{eq17}, we get
\begin{equation}\label{eq18}
\begin{split}
P\{N_\lambda(s)=0\}=1-\sum_{k=1}^\infty P\{N_\lambda(s)=k\}&=1-\sum_{k=1}^\infty \alpha s \frac{(1)_{k,\lambda}}{k!}\theta^k+o(s)\\
& = 1-\alpha s(e_\lambda(\theta)-1)+o(s).
\end{split}
\end{equation}

For $N_\lambda(s)=k \ (k=0,1,2,\cdots)$, we have
\begin{equation}\label{eq19}
\begin{split}
E[\exp(-xN_\lambda(s))]&=\sum_{k=0}^\infty e^{-xk}P\{N_\lambda(s)=k\} \\
&=P\{N_\lambda(s)=0\}+\sum_{k=1}^\infty e^{-xk} P\{N_\lambda(s)=k\} \\
&=1-\alpha s(e_\lambda(\theta)-1)+o(s)+\alpha s \sum_{k=1}^\infty e^{-xk}\frac{(1)_{k,\lambda}}{k!}\theta^k \\
&=1-\alpha s(e_\lambda(\theta)-1)+o(s)+\alpha s (e_\lambda(e^{-x}\theta)-1) \\
&=1+\alpha s(e_\lambda(e^{-x}\theta)-e_\lambda(\theta))+o(s).
\end{split}
\end{equation}

From \eqref{eq16} and \eqref{eq19}, we have
\begin{equation}\label{eq20}
\begin{split}
g(t+s)&=g(t)E[\exp(-xN_\lambda(s))] \\
&=g(t)(1+\alpha s(e_\lambda(e^{-x}\theta)-e_\lambda(\theta)))+o(s).
\end{split}
\end{equation}

Thus, by \eqref{eq20}, we get
\begin{equation}\label{eq21}
\begin{split}
g(t+s)-g(t)=g(t)\big(\alpha s(e_\lambda(e^{-x}\theta)-e_\lambda(\theta)\big)+o(s).
\end{split}
\end{equation}

From \eqref{eq21}, we note that
\begin{equation}\label{eq22}
\begin{split}
\frac{dg(t)}{dt}=\lim_{s\rightarrow0}\frac{g(t+s)-g(t)}{s}=g(t)\alpha(e_\lambda(e^{-x}\theta)-e_\lambda(\theta)).
\end{split}
\end{equation}

Thus, by \eqref{eq22}, we get
\begin{equation}\label{eq23}
\begin{split}
\log{g(t)}=\alpha t(e_{\lambda}(e^{-x}\theta)-e_{\lambda}(\theta)).
\end{split}
\end{equation}

From \eqref{eq23}, we note that
\begin{equation}\label{eq24}
\begin{split}
g(t)=\exp(\alpha t(e_{\lambda}(e^{-x}\theta)-e_{\lambda}(\theta))).
\end{split}
\end{equation}

Thus, by \eqref{eq24}, we get
\begin{equation}\label{eq25}
\begin{split}
E[\exp(-xN_{\lambda}(t))]=g(t)=\exp(\alpha t(e_{\lambda}(e^{-x}\theta)-e_{\lambda}(\theta))).
\end{split}
\end{equation}

Therefore, by Theorem 3 and \eqref{eq25}, we obtain the following theorem.
\begin{theorem}
If $\{N_{\lambda}(t)|\,t\geq0\}$ is a degenerate Bell process with parameters $(\alpha, \theta)$, then
\begin{equation*}
\begin{split}
N_{\lambda}(t) \sim DB_{\lambda}(\alpha t, \theta), \ \ {\rm{for \  all}} \quad t\geq0.
\end{split}
\end{equation*}
\end{theorem}

\medskip

In view of Theorem 5, we may have another definition for the degenerate Bell process.

\begin{definition}
A counting process $\{N_{\lambda}(t)|\,t\geq0\}$ is called a degenerate Bell process with parameters $(\alpha, \theta)$ if the following assumptions hold:
\begin{itemize}
\item[(i)] $N_\lambda(0)=0$, \\
\item[(ii)] $\{N_\lambda(t)|\,t\geq0\}$ has stationary and independent increments, \\
\item[(iii)]for all $k\in\mathbb{N}\cup\{0\}$ and $t>0$,
$P\{N_{\lambda}(t)=k\}=e^{-\alpha t(e_{\lambda}(\theta)-1)}\frac{\theta^k}{k!}\phi_{k,\lambda}(\alpha t)$,
\end{itemize}
\end{definition}
\noindent where $\phi_{k,\lambda}(x)$ are the degenerate Bell polynomials.

\vspace{0.1in}

Let $\{N_{i,\lambda}(t)|\, t\geq 0 \}_{i \in \mathbb{N}}$ be the family of degenerate Bell process with parameters $(\alpha_i,\theta)_{i \in \mathbb{N}}$.
Let $\widetilde{N}_{n,\lambda}(t)=\sum_{i=1}^nN_{i,\lambda}(t)$, for all $t\geq 0$,  and $\beta_n=\sum_{i=1}^n \alpha_i$.

Recalling \eqref{eq07-1}, we obtain
\begin{equation}\label{eq30}
\begin{split}
P[N_{1,\lambda}(t)+N_{2,\lambda}&(t)=k]=\sum_{i=0}^k P[N_{1,\lambda}(t)=i]P[N_{2,\lambda}(t)=k-i] \\
&=\sum_{i=0}^{k} e^{-\alpha_{1}t(e_\lambda(\theta)-1)}\frac{\theta^i}{i!}\phi_{i,\lambda}(\alpha_{1}t)e^{-\alpha_2 t(e_\lambda(\theta)-1)}\frac{\theta^{k-i}}{(k-i)!}\phi_{k-i,\lambda}(\alpha_2 t) \\
&=e^{-(\alpha_1+\alpha_2 )t(e_\lambda(\theta)-1)}\frac{\theta^k}{k!}\sum_{i=0}^k \binom{k}{i}\phi_{i,\lambda}(\alpha_1 t)\phi_{k-i, \lambda}(\alpha_2 t) \\
&=e^{-(\alpha_1+\alpha_2)t(e_\lambda(\theta)-1)}\frac{\theta^k}{k!}\phi_{k,\lambda}((\alpha_1+\alpha_2)t).
\end{split}
\end{equation}
From \eqref{eq30}, we have
\begin{equation}\label{eq31}
\begin{split}
P[\widetilde{N}_{2,\lambda}(t)=k]=e^{-\beta_2 t(e_{\lambda}(\theta)-1)}\frac{\theta^k}{k!}\phi_{k,\lambda}(\beta_2 t).
\end{split}
\end{equation}
For $n \ge 0$, we assume that the probability mass function of $\widetilde{N}_{n,\lambda}(t)$ is given by
\begin{equation}\label{eq29}
\begin{split}
P[\widetilde{N}_{n,\lambda}(t)=k]=e^{-\beta_n t(e_\lambda(\theta)-1)}\frac{\theta^k}{k!}\phi_{k,\lambda}(\beta_n t).
\end{split}
\end{equation}
Then we obtain
\begin{equation}\label{eq32}
\begin{split}
P[\widetilde{N}_{n,\lambda}(t)+&N_{(n+1),\lambda}(t)=k]=\sum_{i=0}^kP[\widetilde{N}_{n,\lambda}(t)=i]P[N_{(n+1),\lambda}(t)=k-i] \\
&=\sum_{i=0}^ke^{-\beta_n t(e_\lambda(\theta)-1)}\frac{\theta^i}{i!}\phi_{i,\lambda}(\beta_n t)e^{-\alpha_{n+1}t(e_\lambda(\theta)-1)}\frac{\theta^{k-i}}{(k-i)!}\phi_{k-i,\lambda}(\alpha_{n+1}t) \\
&=e^{-(\beta_n +\alpha_{n+1})t(e_\lambda(\theta)-1)}\frac{\theta^k}{k!}\phi_{k,\lambda((\beta_n+\alpha_{n+1}t)}.
\end{split}
\end{equation}
By \eqref{eq32}, we get
\begin{equation}\label{eq33}
\begin{split}
P[\widetilde{N}_{n+1}(t)=k]=e^{-\beta_{n+1} t(e_\lambda(\theta)-1)}\frac{\theta^k}{k!}\phi_{k,\lambda}(\beta_{n+1} t).
\end{split}
\end{equation}

From \eqref{eq33}, we have the following result.
\begin{theorem}
Let $\{N_{i,\lambda}(t)|\, t\geq 0\}_{i \in \mathbb{N}}$ be the family of degenerate Bell process with parameters $(\alpha_i,\theta)_{i \in \mathbb{N}}$.
Let $\widetilde{N}_{n,\lambda}(t)=\sum_{i=1}^nN_{i,\lambda}(t)$, for all $t\geq 0$,  and let $\beta_n=\sum_{i=1}^n \alpha_i$. Then
$\{\widetilde{N}_{n,\lambda}(t) | \, t\ge 0\}_{n \in \mathbb{N}}$ is the family of degenerate Bell process with parameters $(\beta_n, \theta)_{n \in \mathbb{N}}$.
\end{theorem}

\begin{remark}
We note that for two Bell process $N_{1}(t)$  and $N_{2}(t)$ with respective parameters $(\alpha_1,\theta_1)$ and $(\alpha_2,\theta_2)$, $N_{1}(t)+N_{2}(t)$ is not a Bell process (see \cite{6}).
\end{remark}

\bigskip

\section{Conclusion}
In recent years, the exploration for degenerate versions has been done for many special numbers and polynomials. It is remarkable that this led to find the degerate gamma functions, the degenerate umbral calculus and the degenerate $q$-umbral calculus. Here the central role is played by the degenerate exponentials in all of these quests (see \eqref{eq01}). \par
In this paper, we considered the discrete degenerate Bell distributions and the degenerate Bell process which are degenerate versions of the Bell-Touchard probability distributions and the Bell-Touchard process, respectively. Several properties were derived for the degenerate Bell distribution. The degenerate Bell process was introduced by giving two equivalent definitions. Then we showed one method of constructing a new infinite family of degenerate Bell process out of a given infinite family of degenerate Bell process. \par
We would like to continue to explore degenerate versions and to find their applications to physics, science and engineering as well as to mathematics.

\vspace{1cm}



\noindent{\bf {Availability of data and material}} \\
Not applicable.

\vspace{0.1in}

\noindent{\bf{Funding}} \\
The third author is supported by the Basic Science Research Program, the National
Research Foundation of Korea, (NRF-2021R1F1A1050151).
\vspace{0.1in}

\noindent{\bf{Ethics approval and consent to participate}} \\
The authors declare that there is no ethical problem in the production of this paper.

%
%
%
%

%

\


\bigskip
\begin{thebibliography}{9}


\bibitem{1} Beichelt, F. Applied probability and stochastic processes. Second edition. CRC Press. Boca Raton, FL, 2016. front matter+560 pp. ISBN: 978-1-4822-5764-960-01

\bibitem{2} Castellares, F.; Lemonte, A. J ; Moreno-Arenas, G. On the two-paramrter Bell-Touchard discrete distribution. Communications in Statistics - Theory and Methods 49 (2020), no. 19, 4834-4852.

\bibitem{3} Comtet, L. Advanced combinatorics. The art of finite and infinite expansions. Revised and enlarged edition. D. Reidel Publishing Co., Dordrecht, 1974. xi+343 pp. ISBN: 90-277-0441-4

\bibitem{4} Djordjevic, G. B.; Milovanovic, G. V.  Special classes of polynomials. University of Nis, Faculty of Technology Leskovac, 2014.

\bibitem{5} Freud, T.; Rodriguez, P. M. The Bell-Tochard counting process. Appl. Math. Comput. 444 (2023), paper no. 127741.

\bibitem{6} J\'{a}nossy, L.; R\'{e}nyi, A.; Acz\'{e}l, J. On composed Poisson distributions, I. Acta Math. Acad. Sci. Hungar. 1 (1950), no. 2-4, 209-224.

\bibitem{7} Kim, H. K. Fully degenerate Bell polynomials associated with degenerate Poisson random variables. Open Math. 19 (2021), no. 1, 284-296.

\bibitem{8} Kim, T.; Kim, D. S. Degenerate zero-truncated Poisson random variables. Russ. J. Math. Phys. 28 (2021), no. 1, 66-72.

\bibitem{9} Kim, T.; Kim, D. S. Some identities on truncated polynomials associated with degenerate Bell polynomials. Russ. J. Math. Phys. 28 (2021), no. 3, 342-355.

\bibitem{10} Kim, T.; Kim, D. S. Some identities on degenerate Bell polynomials and their related identities. Proc. Jangjeon Math. Soc. 25 (2022), no. 1, 1-11.

\bibitem{11} Kim, T.; Kim, D. S.; Dolgy, D. V.; Park, J.-W. Degenerate binomial and Poisson random variables associated with degenerate Lah-Bell polynomials. Open Math. 19 (2021), no. 1, 1588-1597.

\bibitem{12} Ross, S. M. Introduction to probability models. Twelfth edition of [MR0328973]. Academic Press, London, 2019. xv+826 pp. ISBN: 978-0-12-814346-9 60-01

\bibitem{13} Theodorescu, R.; Borwein, J. M. Problems and Solutions: Solutions: Moments of the Poisson distribution: 10738. Amer. Math. Monthly 107 (2000), no. 7, 659.


%





\end{thebibliography}
\end{document}